\documentclass[a4paper, 11pt]{article}
\usepackage[utf8]{inputenc}
\usepackage[margin=3cm]{geometry} 

\usepackage{amsmath,amsthm,amssymb}
\usepackage{enumerate}
\usepackage{verbatim}
\usepackage{graphicx}
\usepackage{listings}
\usepackage{longtable}
\usepackage{array}
\usepackage{float}
\usepackage{caption}
\usepackage{subcaption}
\usepackage{bm}
\usepackage{amsthm}
\usepackage{afterpage}
\usepackage{tikz}
\usepackage{tikz-cd}
\usepackage{standalone}
\usepackage{booktabs}

\allowdisplaybreaks  
\newcounter{mainnum}

\title{Potential good reduction of hyperelliptic curves}
\author{Robin Visser}
\date{19 January 2022}

\begin{document}

\maketitle

\begin{abstract}
    Let $K$ be a number field, and $g \geq 2$ a positive integer. We define $c_K(g)$ as the smallest integer $n$ such that there exist infinitely many $\overline{K}$-isomorphism classes of genus $g$ hyperelliptic curves $C/K$ with all Weierstrass points in $K$ having potentially good reduction outside $n$ primes in $K$.  We show that $c_K(g) > \pi_{K, \textrm{odd}}(2g) + 1$, where $\pi_{K, \textrm{odd}}(n)$ denotes the number of odd primes in $K$ with norm no greater than $n$, as well as present a summary of various conditional and unconditional results on upper bounds for $c_K(g)$.

\end{abstract}

{\renewcommand{\thefootnote}{}
}

\section{Introduction}

For a given number field $K$, the problem of studying the reduction behaviour of hyperelliptic curves $C/K$ has attracted significant interest over the last few decades. Indeed, there exist well-known effective algorithms for determining all hyperelliptic curves $C/K$ of a given genus $g$ with potentially good reduction outside a given finite set of primes $S$, with Merriman--Smart \cite{merrimansmart, smartbinaryforms} effectively carrying out this computation to list all genus 2 curves $C/\mathbb{Q}$ with good reduction outside 2. \\

For our purposes, we shall be interested in the reduction of hyperelliptic curves $C/K$ with all of its Weierstrass points lying in $K$. We first recall that a prime $\mathfrak{p}$ in $K$ is considered \textrm{odd} if it lies above an odd rational prime (or equivalently has odd norm), and define $\pi_{K, \textrm{odd}}(n)$ as the number of odd primes in $K$ with norm no greater than $n$. We also define $\mathcal{B}_\textrm{odd}(C/K)$ as the set of odd primes $\mathfrak{p}$ in $K$ for which $C/K$ does not have potential good reduction at $\mathfrak{p}$. Our main result is the following: \\

\refstepcounter{mainnum} \label{thm:maintheorem}
\textbf{Theorem \themainnum.} (\textit{Theorem \ref{thm:norm}, Theorem \ref{thm:mainlowerboundtheorem}}) Let $K$ be a number field.  Then for all genus $g$ hyperelliptic curves $C/K$ with all Weierstrass points in $K$, we have $\mathfrak{p} \in \mathcal{B}_\textrm{odd}(C/K)$ for all odd primes $\mathfrak{p}$ satisfying $N_{K/\mathbb{Q}}(\mathfrak{p}) < 2g$.
Furthermore, there are only finitely many $\overline{K}$-isomorphism classes of genus $g$ hyperelliptic curves $C/K$ with all Weierstrass points in $K$ satisfying $\# \mathcal{B}_\textrm{odd}(C/K) \leq \pi_{K, \textrm{odd}}(2g) + 1$. \\

This gives us the lower bound $c_K(g) > \pi_{K, \textrm{odd}}(2g) + 1$. Applying Theorem \ref{thm:maintheorem} to $K = \mathbb{Q}$ gives the following corollaries: \\

\refstepcounter{mainnum}
\textbf{Corollary \themainnum.} (\textit{Corollary \ref{thm:rationalWeier}, Corollary \ref{thm:box}}) Let $C/\mathbb{Q}$ be a hyperelliptic curve with rational Weierstrass points.  Then $C$ cannot have potentially good reduction at any odd prime $p \leq 2g$. Moreover, there is no genus 2 curve $C/\mathbb{Q}$ with rational Weierstrass points having potentially good reduction outside one prime. \\

Finally, in Section 4, we also prove the following various conditional and unconditional upper bounds for $c_K(g)$. \\

\refstepcounter{mainnum}
\label{thm:mainupperboundtheorem}
\textbf{Theorem \themainnum.} (\textit{Theorem \ref{thm:upperDickson}, Theorem \ref{thm:upper}, Theorem \ref{thm:uppergenDickson}}) Let $K$ be a number field of degree $n$. Then $c_K(g) \leq (\frac{2}{\log{2}} + o(1)) n g \log{g}$.  Furthermore, under the assumption of the Hardy-Littlewood prime $k$-tuples conjecture for $K$, we have $c_K(g) \leq 2g - 1 + n \pi(2g)$, and under the assumption of Schinzel’s hypothesis H for $K$, we have moreover that
\begin{equation*}
    c_K(g) \leq \sum_{\substack{1 \leq d < g, \textrm{ or} \\ d < 2g, \, d \textrm{ even}}} \frac{n}{[K(\zeta_d) : \mathbb{Q}(\zeta_d)]} + 1 +  n \pi(2g) .
\end{equation*}

Precise statements of the Hardy--Littlewood prime $k$-tuples conjecture and the Schinzel hypothesis H are provided in Section 4. This hence gives rise to the following corollaries: \\

\refstepcounter{mainnum}
\textbf{Corollary \themainnum.} (\textit{Corollary \ref{thm:upperCor1}, Corollary \ref{thm:upperCor2}}) Let $K$ be a number field of degree $n$ with no non-trivial abelian subfields, and suppose that Schinzel's hypothesis H holds for $K$. Then if $K$ is abelian (and hence of prime degree) with conductor $\mathfrak{f}_K$, then
\begin{equation*}
    c_K(g) \leq \begin{cases}
    \frac{3}{2} g \big(1 + \frac{n-1}{\mathfrak{f}_K} \big) + 1 + n \pi(2g) & \textrm{ if $\mathfrak{f}_K$ odd, } \\[2mm]
    \frac{3}{2} g \big( 1 + \frac{4(n-1)}{3 \mathfrak{f}_K} \big) + 1 + n \pi(2g) & \textrm{ if $\mathfrak{f}_K$ even, } 
    \end{cases}
\end{equation*}
otherwise $c_K(g) \leq \frac{3}{2} g + n \pi (2g)$ if $K$ is non-abelian. \\

Let $C/K$ be a hyperelliptic curve of genus $g$ with all Weierstrass points in $K$.  We first recall that $C$ is isomorphic to a curve of the form
\begin{equation*}
    y^2 = cx(x-1)(x-a_1)(x-a_2) \cdots (x-a_{2g-1})
\end{equation*}
for some $c, a_i \in K$, $a_i \not = 0, 1$, as all Weierstrass points of $C$ lie in $K$. This is called the \textit{Rosenhain normal form} of $C$. \\

We shall now apply the machinery of cluster pictures to study these curves, introduced by Dokchitser, Dokchitser, Maistret, Morgan \cite{dokchitser}.

\section{Cluster pictures}

We first quickly recall the notion of cluster pictures.  Indeed, let $C/K$ be a hyperelliptic curve, given by $C : y^2 = f(x)$, and let $\mathcal{R}$ denote the roots of $f$, with $\mathcal{P}(\mathcal{R})$ being the power set of $\mathcal{R}$. Let $\mathfrak{p}$ be an odd prime in $K$, and let $v_\mathfrak{p}$ denote the discrete normalised $p$-adic valuation induced by $\mathfrak{p}$. The \textbf{cluster picture} $\Sigma_{\mathfrak{p}} \subset \mathcal{P}(\mathcal{R})$ associated to $C$ (with respect to $\mathfrak{p}$) is given by the following set:
\begin{equation*}
    \Sigma_{\mathfrak{p}} :=  \big\{ \mathfrak{s} \in \mathcal{P}(\mathcal{R}) \;|\; \forall x \in \mathfrak{s}, \, v_{\mathfrak{p}}(x - z) \geq d \textrm{ for some } z \in \overline{K}, d \in \mathbb{Q} \big\} .
\end{equation*}

We say that $\Sigma_{\mathfrak{p}}$ is \textit{trivial} if it only contains the singleton elements and $\mathcal{R}$ itself. \\

Using these cluster pictures, we shall make use of the following theorem of Dokchitser, Dokchitser, Maistret, Morgan, which allows us to directly read off whether a hyperelliptic curve $C/K$ has potentially good reduction at $\mathfrak{p}$. \\

\refstepcounter{mainnum} \label{thm:potgoodredn}
\textbf{Theorem \themainnum.} \cite[p.~4]{dokchitser2} Let $C/K$ be a hyperelliptic curve of genus $g$, and let $\mathfrak{p}$ be an odd prime in $K$. Then $C$ has potentially good reduction at $\mathfrak{p}$ if and only if $\Sigma_{\mathfrak{p}}$ has no proper clusters of size $< 2g + 1$. \\

From this theorem, we can prove the following proposition. \\

\refstepcounter{mainnum} \label{thm:lambdas}
\textbf{Proposition \themainnum.} Let $C/K$ be a hyperelliptic curve with Weierstrass points in $K$, given in Rosenhain normal form
\begin{equation*}
    y^2 = c x(x-1)(x-\lambda_1) \dots (x-\lambda_{2g-1}) , \qquad c, \lambda_i \in K .
\end{equation*}
Let $\mathfrak{p}$ be an odd prime of $K$. Then $C$ has potentially good reduction at $\mathfrak{p}$ if and only if we have $v_\mathfrak{p}(\lambda_i) = v_\mathfrak{p}(\lambda_i - 1) = 0$ for all $i \in 1, \dots, 2g-1$, and $v_\mathfrak{p}(\lambda_i - \lambda_j) = 0$ for all distinct $i, j \in \{1, \dots, 2g-1\}$. (i.e. the values $\lambda_i, \lambda_i - 1, \lambda_i - \lambda_j$ are all $\mathfrak{p}$-units) \\

\textit{Proof.} Let $C/K$ be given in the above form, and let $\mathcal{R}$ denote the Weierstrass points, i.e. $\mathcal{R} := \{0, 1, \lambda_1, \dots, \lambda_{2g-1} \} $.  Then by Theorem \ref{thm:potgoodredn}, since $|\mathcal{R}| = 2g+1$, we have that $C$ has potentially good reduction at $\mathfrak{p}$ if and only if $\Sigma_\mathfrak{p}$ is trivial. \\

Note that $\Sigma_\mathfrak{p}$ is trivial if and only if $v_\mathfrak{p}(r_i - r_j)$ is constant over all distinct pairs $r_i, r_j \in \mathcal{R}$.  However, since $v_\mathfrak{p}(1 - 0) = 0$, this implies that $v_\mathfrak{p}(\lambda_i) = v_\mathfrak{p}(\lambda_i - 1) = 0$ for all $i$, and that $v_\mathfrak{p}(\lambda_i - \lambda_j) = 0 $ for all $i, j$, which yields the result. \qed \\

This immediately implies the following corollary: \\

\refstepcounter{mainnum} \label{cor:Cpgr}
\textbf{Corollary \themainnum.} Let $K$ be a number field, and let $S$ be a finite set of primes of $K$, and assume that $S$ consists of all even primes of $K$. Let $\mathcal{O}_S^\times$ denote the set of $S$-units in $K$.  Then for a given hyperelliptic curve $C/K$ of the above form, $C$ has potentially good reduction outside $S$ if and only if $\lambda_i$ and $\lambda_i - 1$ are in $\mathcal{O}_S^\times$, and if $\lambda_i - \lambda_j$ are in $\mathcal{O}_S^\times$. \\

We can now prove our main theorem. \\

\refstepcounter{mainnum} \label{thm:norm}
\textbf{Theorem \themainnum.} Let $C/K$ be a hyperelliptic curve with Weierstrass points in $K$. Then $C$ cannot have potentially good reduction at any odd prime $\mathfrak{p}$ such that $N_{K/\mathbb{Q}}(\mathfrak{p}) \leq 2g$. \\

\textit{Proof.} Let $C$ be given by its Rosenhain normal form:
\begin{equation*}
    y^2 = cx(x-1)(x-\lambda_1) \cdots (x-\lambda_{2g-1})
\end{equation*}
and assume for contradiction that $C$ has potentially good reduction at $\mathfrak{p}$, where $\mathfrak{p}$ is a prime ideal of $K$ such that $N(\mathfrak{p}) \leq 2g$. \\

We have by Theorem \ref{thm:lambdas} that $\lambda_1, \dots, \lambda_{2g-1}$ must all be $\mathfrak{p}$-units.  Furthermore, note that each of the roots $0, 1, \lambda_1, \dots, \lambda_{2g-1}$ must yield distinct values under the reduction map $\mathcal{O}_K \to \mathcal{O}_K / \mathfrak{p}$. However, this is a contradiction if $2g + 1 > N(\mathfrak{p})$, noting that $\# \mathcal{O}_K / \mathfrak{p} = N(\mathfrak{p})$. \qed \\

We remark that the above inequality is tight, since given any prime $\mathfrak{p}$ with $N(\mathfrak{p}) > 2g$, we can simply let $0, 1, \lambda_1, \dots, \lambda_{2g-1}$ be some distinct representative elements in the residue field to yield an example of a curve $C$ with good reduction at $\mathfrak{p}$. \\

This result immediately implies that there are only finitely many $\overline{K}$-isomorphism classes of genus $g$ hyperelliptic curves $C/K$ with Weierstrass points in $K$ having potentially good reduction outside at most $\pi_{K, \textrm{odd}}(2g)$ odd primes. However, remarkably we can go one step further: \\

\refstepcounter{mainnum} \label{thm:mainlowerboundtheorem}
\textbf{Theorem \themainnum.} There are only finitely many $\overline{K}$-isomorphism classes of genus $g$ hyperelliptic curves $C/K$ with rational Weierstrass points in $K$ having potentially good reduction outside at most $\pi_{K, \textrm{odd}}(2g) + 1$ primes. \\

In order to prove the above theorem, we shall make use of the following elementary (albeit technical) lemma: \\

\refstepcounter{mainnum} \label{thm:mainlemma}
\textbf{Lemma \themainnum.} Let $K$ be a number field, and $S$ a fixed finite set of primes of $K$. Then there exist only finitely many odd primes $\mathfrak{p}$ such that there exist distinct $T$-units $x, y, z \in \mathcal{O}_T^\times$ where $T = S \cup \{\mathfrak{p} \}$ such that $x-y$, $x-z$, and $y-z$ are all $T$-units, and assuming $v_\mathfrak{p}(x), v_\mathfrak{p}(y), v_\mathfrak{p}(z), v_\mathfrak{p}(x-y), v_\mathfrak{p}(x-z), v_\mathfrak{p}(y-z)$ are not all equal. \\

The proof of this lemma proceeds by analysing various three-term $S$-unit equations, and thus we shall make essential use of the following finiteness result: \\

\refstepcounter{mainnum} \label{thm:sunit3finite}
\textbf{Theorem \themainnum.} \cite[p.~131]{sunit3} Let $K$ be a number field and $S$ a fixed finite set of primes of $K$. Then the equation $u+v+w = 1$ has only finitely many solutions in $u, v, w \in \mathcal{O}_S^*$, such that $u, v, w \not = 1$ (i.e. only finitely many non-degenerate solutions). \\

With the above theorem under our belt, we can now prove Lemma \ref{thm:mainlemma}. \\

\textit{Proof of Lemma \ref{thm:mainlemma}.} We first fix some odd prime $\mathfrak{p}$ of $K$, and shall aim to derive a set of equations which can only be satisfied for finitely many $\mathfrak{p}$. With this in mind, we can first assume without loss of generality that $v_\mathfrak{p}(x) \geq v_\mathfrak{p}(y) \geq v_\mathfrak{p}(z)$.  For brevity we denote $s := \frac{x}{z}$ and $t := \frac{y}{z}$ and define $a := v_\mathfrak{p}(s)$ and $b := v_\mathfrak{p}(t)$, noting that $a, b$ are non-negative integers where $a \geq b$. 
We also define $c := v_\mathfrak{p}(s - 1)$, $d := v_\mathfrak{p}(t-1)$ and $e := v_\mathfrak{p}(s-t)$, noting that $s, t$ and $s-t$ are all $T$-units. \\

The proof now proceeds by considering the various cases for $a$ and $b$.  In each case, the main idea is to obtain a three term $S$-unit equation from which we will obtain only finitely many solutions.

\begin{itemize}
    \item \textbf{Case 1: } $a, b > 0$ and $a > b$. This implies $c = d = 0$ and $e = b$, and thus we have
    \begin{equation*}
        s - 1 = u, \quad \textrm{and} \quad 
        t - 1 = v,
    \end{equation*}
    for some $u, v \in \mathcal{O}_S^\times$. As $v_\mathfrak{p}(u-v) = v_\mathfrak{p}(t)$, we have that $t/(u-v)$ is an $S$-unit, and thus by rearranging, we obtain the three term $S$-unit equation:
    \begin{equation} \label{eq:case1}
        \frac{t}{u-v} u - \frac{t}{u-v} v - v = 1 .
    \end{equation}
    
    At this stage, we would like to apply Theorem \ref{thm:sunit3finite} in order to conclude that there are only finitely many solutions to the above equation. We must therefore check that we do not obtain (or only obtain finitely many) degenerate solutions where one of the above terms equals 1:

    \begin{itemize}
        \item[(i)] If the first term of (\ref{eq:case1}) is 1, then  $t = v-u$ which implies $x = 0$, contradiction.
        
        \item[(ii)] If the second term of (\ref{eq:case1}) is 1, then $t v = v-u$, which implies $s = t (1-v)$ and hence $1-v$ has positive $\mathfrak{p}$-adic valuation. But $v-1 = t - 2$ which yields a contradiction, since $\mathfrak{p}$ is odd.
        
        \item[(iii)] If the third term of (\ref{eq:case1}) is 1, then $t = 0$, contradiction.
    \end{itemize}
    
    Therefore, we have a three-term $S$-unit equation with no degenerate solutions. Thus, there are only finitely many solutions to (\ref{eq:case1}), and thus only finitely many $v$, and thus clearly only finitely many $\mathfrak{p}$, noting that $b$ is positive. 

    \item \textbf{Case 2: } $a, b > 0$ and $a = b$. As before, we have
    \begin{equation*}
        s - 1 = u, \quad \textrm{and} \quad 
        t - 1 = v
    \end{equation*}
    for some $u, v \in \mathcal{O}_S^\times$. Noting that $v_\mathfrak{p}(t/s) = 0$, by rearranging, we obtain the three term $S$-unit equation:
    \begin{equation} \label{eq:case2}
        \frac{t}{s} u + \frac{t}{s} - v = 1
    \end{equation}
    
    Once again, we check the three cases: 
    \begin{itemize}
        \item[(i)] If the first term of (\ref{eq:case2}) is 1, then $tu = s$ and $t = vs$ which implies $uv = 1$. Now by multiplying the first two equations we get
        \begin{equation*}
            st - (s+t) + 1 = (s-1)(t-1) = uv = 1
        \end{equation*}
        which implies $v_\mathfrak{p}(s+t) = v_\mathfrak{p}(st) = 2a > a$, and thus $v_\mathfrak{p}(s-t) = v_\mathfrak{p}(s+t - 2t) = a$ as $\mathfrak{p}$ odd. This thus yields the following two-term $S$-unit equation:
        \begin{equation*}
            \frac{s}{t} - \frac{s-t}{t} = 1
        \end{equation*}
        which implies finitely many values for $\frac{s}{t}$ and thus for $u$, and so only finitely many values for $\mathfrak{p}$.
        
        \item[(ii)] If the second term of (\ref{eq:case2}) is 1, then $x = y$, contradiction.
        
        \item[(iii)] If the third term of (\ref{eq:case2}) is 1, then $u = 1$ and thus $x/z = 2$, contradiction.
        
    \end{itemize}
    
    Therefore, as before, only finitely many solutions.
    
     \item \textbf{Case 3: } $a > 0$ and $b = 0$. We therefore have $c = 0$ and $e = 0$. This yields
     \begin{equation*}
         s - 1 = u, \quad
         s - t = w
     \end{equation*}
     for some $u, w \in \mathcal{O}_S^\times$, which yields the three term $S$-unit equation:
     \begin{equation} \label{eq:case3}
         w + t - u = 1
     \end{equation}
     
     \begin{itemize}
        \item[(i)] If $w = 1$, then $v_\mathfrak{p}(t-1) = v_\mathfrak{p}(s - 2) = 0$, as $\mathfrak{p}$ odd. Therefore $t - 1 = v$ for some $v \in \mathcal{O}_S^\times$. This yields a 2-term $S$-unit equation, of which there are only finitely many solutions for $t, v$, and thus for $u$, hence only finitely many for $p$.
        
        \item[(ii)] If $t = 1$, then $y = z$, contradiction.
        
        \item[(iii)] If $u = -1$. then $x = 0$, contradiction.

    \end{itemize}
    
    \item \textbf{Case 4: } $a = b = 0$ and $c > d$. This implies $e = d$ which yields the three term $S$-unit equation:
    \begin{equation} \label{eq:case4}
        \frac{t-1}{s-t} t - \frac{t-1}{s-t} s + t = 1
    \end{equation}
    
    Firstly, if $d = 0$, then $v_\mathfrak{p}(t-1) = 0$ which yields a two term $S$-unit equation of which there are only finitely many solutions. Thus, we may assume $d > 0$.
    
    \begin{itemize}
        \item[(i)]  If the first of (\ref{eq:case4}) is 1, then $s - t = t(t-1)$ which implies $s = t^2$. This yields
        \begin{equation*}
            s - 1 = (t-1)(t+1)
        \end{equation*}
        which implies $t+1$ has positive $\mathfrak{p}$-adic valuation. But $t + 1 = (t-1) + 2$ which yields a contradiction as $\mathfrak{p}$ odd.
        
        \item[(ii)] If the second term of (\ref{eq:case4}) is 1, then $t-1 = t - s$ and so $s = 1$, contradiction.
        
        \item[(iii)] If the third term of (\ref{eq:case4}) is 1, then $y = z$, contradiction.
        
    \end{itemize}
    
    \item \textbf{Case 5: } $a = b = 0$ and $c = d$. Again,  note that if $c = d = 0$, then $s$ and $t$ satisfy two-term $S$-unit equations, of which there are only finitely many solutions. This thus implies only finitely many $p$, since we'd then have $v_\mathfrak{p}(s-t) > 0$ by assumption.
    
    Now assume $c, d \not = 0$, and note that $c, e$ must necessarily be positive. We obtian the three term $S$-unit equation:
    \begin{equation} \label{eq:case5}
        \frac{s-1}{t-1} - \frac{s-1}{t-1}t + s = 1
    \end{equation}
    
    \begin{itemize}
        \item[(i)] If the first term of (\ref{eq:case5}) is 1, then $s = t$, contradiction.
        
        \item[(ii)] If the second term of (\ref{eq:case5}) is 1, then $(s-1)t = - (t-1)$ and $(s-1) = -s(t-1)$ which implies $st = 1$.  By a dual argument to the above case 4(i), we have $v_\mathfrak{p}(t+1) = 0$. This implies
        \begin{equation*}
            e = v_\mathfrak{p}(s - t) = v_\mathfrak{p}(1-t^2) = v_\mathfrak{p}(1-t) = c
        \end{equation*}
        This implies that we have the following two term $S$-unit equation:
        \begin{equation*}
            \frac{s-1}{t-1} - \frac{s-t}{t-1} = 1
        \end{equation*}
        which implies only finitely many solutions for $\frac{s-1}{t-1}$.  Therefore, this gives finitely many $t$, and thus finitely many $\mathfrak{p}$.
        
        \item[(iii)] If $s = 1$, then $x = z$, contradiction.
        
    \end{itemize}
    
    \item \textbf{Case 6: } $a = b = 0$ and $c < d$. Done analogously to case 4.

\end{itemize}

Therefore, in each case, we obtain only finitely many valid primes $\mathfrak{p}$, which concludes the proof. \qed \\

We note that effectively obtaining a list of all possible primes $\mathfrak{p}$ depends entirely on the effectiveness of solving the above three term $S$-unit equations. From the results of \cite{sunit3}, no finite algorithm has been found to determine all possible solutions, however one can obtain an explicit bound on the number of possible $\mathfrak{p}$, which for a fixed number of terms, is exponential in $|S|$ \cite[p.~132]{sunit3}. \\

We are now finally ready to prove Theorem \ref{thm:mainlowerboundtheorem}: \\

\textit{Proof of Theorem \ref{thm:mainlowerboundtheorem}.}  Let $C/K$ be a hyperelliptic curve of genus $g$ given in Rosenhain normal form $C : y^2 = x(x-1)(x-\lambda_1) \cdots (x - \lambda_{2g-1})$
with Weierstrass points in $K$. By Theorem \ref{thm:norm}, $C$ cannot have potentially good reduction at any odd primes with norm less than $2g$.  Now assume $C$ has potentially good reduction outside exactly $\pi_{K,\textrm{odd}}(2g)+1$ odd primes $S$. \\

Thus, $S$ must consist of all $\pi_{K,\textrm{odd}}(2g)$ primes with norm below $2g$, plus one additional prime $\mathfrak{p}$.  Now by Corollary \ref{cor:Cpgr}, we must have that $\lambda_1, \lambda_2, \lambda_1 - 1, \lambda_2 - 1$ and $\lambda_1 - \lambda_2$ are all $S$-units.  Therefore, by Lemma \ref{thm:mainlemma}, there are only finitely many possible primes $\mathfrak{p}$, and thus by either applying a theorem of Faltings \cite[p.~25]{cornellsilverman} or by finiteness of solutions to $S$-unit equations, we obtain only finitely many $\overline{K}$-isomorphism classes of hyperelliptic curves with potentially good reduction outside $S$. \qed \\

We do remark that only the Weierstrass points $0, 1, \lambda_1,$ and $\lambda_2$ were used in the proof above, whilst Corollary \ref{cor:Cpgr} does include constraints on all the Weierstrass points $\lambda_i$.  Indeed, if we were to use all $\lambda_i$, we would expect to prove a significantly stronger lower bound for $c_K(g)$. 

A heuristic argument suggests that if we generalise Lemma \ref{thm:mainlemma} where we adjoin an additional $k$ primes $T := S \cup \{ \mathfrak{p}_1, \dots, \mathfrak{p}_k \}$, then assuming we don't encounter any degenerate solutions, this could yield a potential linear lower bound of $g + \pi_{K, \textrm{odd}}(2g)$ for $c_K(g)$. However, besides extending the case bash analysis for small values of $k$, we do not know at this stage how to produce such a proof for arbitrary $k$. \\

We shall now restrict our attention to the specific case where $C$ is a hyperelliptic curve over $\mathbb{Q}$, with all of its Weierstrass points in $\mathbb{Q}$.

\section{Hyperelliptic curves with rational Weierstrass points}

Firstly, it's worth stating the application of Theorem \ref{thm:norm} to the rational case: \\

\refstepcounter{mainnum} \label{thm:rationalWeier}
\textbf{Corollary \themainnum.} Let $C/\mathbb{Q}$ be a hyperelliptic curve with rational Weierstrass points.  Then $C$ cannot have potentially good reduction at any odd prime $p \leq 2g$. \\

Note that this clearly implies that no genus 2 hyperelliptic curve with rational Weierstrass points has potentially good reduction at 3. This corollary can be applied to give a short proof of the following result from Box and Le Fourn \cite{box}, which was originally proven using a two-dimensional analogue of Baker's and Runge's method applied to the Siegel variety $A_2(2)$ (i.e. the moduli space of principally polarised abelian surfaces with full $2$-torsion).\\

\refstepcounter{mainnum} \label{thm:box}
\textbf{Corollary \themainnum.} \cite[p.~3]{box} There is no genus 2 hyperelliptic curve $C$ over $\mathbb{Q}$ such that all Weierstrass points of $C$ are rational and $C$ has potentially good reduction at all but one of the primes. \\

\textit{Proof.} As shown above, such a curve $C$ cannot have potentially good reduction at $3$.  Now assume for contradiction such a curve has potentially good reduction outside $3$.  By applying Corollary \ref{cor:Cpgr}, we can now effectively compute all genus 2 curves $C/\mathbb{Q}$ with rational Weierstrass points having potentially good reduction outside $S = \{2, 3\}$. \\

By Corollary \ref{cor:Cpgr}, we proceed by solving the $S$-unit equation $x+y = 1$, where $x, y \in \mathcal{O}_S^\times$.
These solutions can be computed using existing algorithms, such as those described by von K\"{a}nel and Matschke \cite{vonkanel}. Using their Sage \cite{sage} implementation, we obtained 21 solutions, and can conclude that any such curve must be isomorphic to one of the following curves:
\begin{align*}
    C_1 : y^2 &= cx(x-1)(x-2)(x-3)(x-4), \quad \textrm{with } \Delta_\textrm{min} = 2^{18} 3^4 \\
    C_2 : y^2 &= cx(x-2)(x-3)(x-4)(x-6), \quad \textrm{with } \Delta_\textrm{min} = 2^{14} 3^6
\end{align*}

where $c \in \mathbb{Z}$ is some squarefree integer. \\

After verifying that each of the curves above do not have potentially good reduction at 2, this gives us our contradiction, and thus the result holds. \qed \\

\section{Upper bounds}

It's worth mentioning some of the results we can obtain regarding upper bounds for $c_K(g)$. Most of the more interesting results are conditional on various conjectures concerning the distribution of primes. \\ 

For a given number field $K$, we recall that a $k$-tuple $(h_1, \dots, h_k)$ of distinct elements in $\mathcal{O}_K$ is \textit{admissible} if the set $\{h_1, \dots, h_k\}$ does not consist of all residues mod $\mathfrak{p}$, for every prime $\mathfrak{p}$ in $K$. We also say that an element $x \in \mathcal{O}_K$ is \textit{prime} if the principal ideal generated by $x$ is prime. We now first recall the Hardy-Littlewood prime $k$-tuples conjecture for $K$: \\

\textbf{Conjecture.} (\textit{Hardy-Littlewood prime $k$-tuples conjecture for number fields}) Let $K$ be a number field and $(h_1, \dots, h_k)$ an admissible $k$-tuple in $\mathcal{O}_K$. Then there exist infinitely many $x \in \mathcal{O}_K$ such that each of $x+h_1, \dots, x+h_k$ is prime. \\

Notably, one can therefore prove the following result, conditional on the assumption of the above conjecture. \\

\refstepcounter{mainnum} \label{thm:upperDickson}
\textbf{Theorem \themainnum.} Let $K$ be a number field of degree $n$. Under the assumption of the Hardy-Littlewood prime $k$-tuples conjecture for $K$, then $c_K(g) \leq 2g - 1 + n \pi(2g)$.  \\

\textit{Proof.} For a given genus $g \geq 2$, we consider the following admissible prime $2g-1$ tuple $(h_1, \dots, h_{2g-1})$:
\begin{equation} \label{eq:admissibletuple}
    (0, \; (2g)!, \; 2 \cdot (2g)!, \; 3 \cdot (2g)!, \; \dots, (2g-2) \cdot (2g)! )
\end{equation}

Now by the prime $k$-tuples conjecture, there exist infinitely many primes $p$ in $K$ such that $p+h_1, \dots, p+h_{2g-1}$ are all prime. Thus, for each such prime $p$, we can construct the genus $g$ hyperelliptic curve $C_p/K$ as
\begin{equation*}
    C_p : y^2 = x(x-p-h_1)(x-p-h_2) \cdots (x-p-h_{2g-1}) (x - 2p - 2g \cdot (2g)! ) .
\end{equation*}

As the only possible primes of bad reduction are those which divide the differences between Weierstrass points, it's clear that the only possible primes of bad reduction are either the primes $p+h_1, \dots, p+h_{2g-1}$ or the primes dividing $(2g)!$, of which there are at most $n \pi(2g)$. \\

This therefore yields the existence of infinitely many genus $g$ hyperelliptic curves $C/K$ satisfying $\# \mathcal{B}_\textrm{odd} (C/K) \leq 2g-1 + n \pi(2g)$. Furthermore, as this yields infinitely many different sets of bad primes $\mathcal{B}_\textrm{odd} (C/K)$, this gives rise to infinitely many $\overline{K}$-isomorphism classes of such curves. This therefore yields the conditional bound $c_K(g) \leq 2g-1 + n \pi(2g)$. \qed \\

Whilst the above result gives a conditional linear upper bound for $c_K(g)$, it's worth noting that we can also give an unconditional linearithmic bound for $c_K(g)$. \\

\refstepcounter{mainnum} \label{thm:upper}
\textbf{Theorem \themainnum.} Let $K$ be a number field of degree $n$. We have $c_K(g) \leq (\frac{2}{\log{2}} + o(1)) n g \log{g}$. \\

\textit{Proof.} Let $(h_1, \dots, h_{2g-1})$ be the same admissible prime tuple as given in (\ref{eq:admissibletuple}).  We shall apply the result of Murty and Vatwani \cite[p.~183]{murtyvatwani}, which asserts the existence of infinitely many integers $k$ such that $(k + h_1) \cdots (k + h_{2g-1})$ has at most $(\frac{2}{\log{2}} + o(1)) g \log{g}$ prime divisors in $\mathbb{Q}$.  By therefore considering the set of genus $g$ hyperelliptic curves
\begin{equation*}
    C_k : y^2 = x(x-k-h_1) \cdots (x-k-h_{2g}) (x - 2k - 2g \cdot (2g)!)
\end{equation*}
this yields the desired upper bound. \qed \\

It's tempting to ask how far we can push our conditional upper bounds.  Whilst a sublinear bound is almost certainly out of reach, we can sharpen the above theorem if we furthermore assume the following generalisation to the Hardy-Littlewood prime $k$-tuples conjecture. This goes by various different names, often called Schinzel's hypothesis H, generalised Dickson's conjecture, or the generalised Bunyakovsky conjecture. \\

\textbf{Conjecture.} (\textit{Schinzel's hypothesis H for number fields}) Let $K$ be a number field and $(f_1, \dots, f_k)$ a collection of $k$ distinct nonconstant irreducible polynomials in $\mathcal{O}_K[x]$, such that 
for all primes $\mathfrak{p}$ in $K$, there exists an $n \in \mathcal{O}_K$ where $v_\mathfrak{p}(f_1(n) f_2(n) \cdots f_k(n)) = 0$.
There there exist infinitely many $x \in \mathcal{O}_K$ such that each of $f_1(x), \dots, f_k(x)$ is prime. \\

Under the assumption of the above conjecture, we can prove the following sharpened upper bound for $c_K(g)$. \\

\refstepcounter{mainnum} \label{thm:uppergenDickson}
\textbf{Theorem \themainnum.} Let $K$ be a number field of degree $n$. Assuming Schinzel's hypothesis H for $K$, we have that
\begin{equation} \label{eq:uppermaingenbound}
    c_K(g) \leq \sum_{\substack{1 \leq d < g, \textrm{ or} \\ d < 2g, \, d \textrm{ even}}} \frac{n}{[K(\zeta_d) : \mathbb{Q}(\zeta_d)]} + 1 +  n \pi(2g)
\end{equation}

\textit{Proof.} For brevity, we shall denote $\alpha := (2g)!$. The idea is to consider, for infinitely many $k$, genus $g$ hyperelliptic curves of the form
\begin{equation*}
    C_k : y^2 = x(x-1)(x+1)(x-\alpha k)(x+\alpha k)(x-(\alpha k)^2)(x+(\alpha k)^2) \cdots (x- (\alpha k)^{g-1} ) (x + (\alpha k)^{g-1} )
\end{equation*}
We note that the only possible primes of bad reduction are those which divide $2 \alpha k$, $(\alpha k)^d - 1$, or $(\alpha k)^d + 1$ for some $d < g$.  Under the assumption of Schinzel's hypothesis H, it thus suffices to count the number of irreducible factors of $(\alpha x)^d \pm 1$ over $K$.

We know that $(\alpha x)^d \pm 1$ factorises over $\mathbb{Q}$ as
\begin{equation*}
    (\alpha x)^d - 1 = \prod_{i | d} \Phi_i(\alpha x), \quad \textrm{and} \quad (\alpha x)^d + 1 = \prod_{\substack{i | 2d \\ i \not | d}} \Phi_i(\alpha x),
\end{equation*}
where $\Phi_i(x)$ denotes the $i$-th cyclotomic polynomial. Furthermore, the factorisation of $\Phi_i(x)$ over $K$ can be given as $\Phi_i(\alpha x) = f_{i,1}(\alpha x) \cdots f_{i, \ell_i}(\alpha x)$ where each $f_{i,j}(x)$ has degree $[K(\zeta_i) : K]$, and hence $\ell_i = \frac{\varphi(i)}{[K(\zeta_i) : K]} = \frac{n}{[K(\zeta_d) : \mathbb{Q}(\zeta_d)]}$ by tower law. \\

Now Schinzel's hypothesis H states that we can find infinitely many primes $p$ in $K$ such that $f_{i,j}(\alpha p)$ are all prime, noting that the factor of $\alpha$ ensures we have no local obstructions to primality. By thus counting the primes dividing $(\alpha p)^d \pm 1$, the prime $p$, and the primes dividing $2 \alpha$, this therefore yields the conditional bound
\begin{equation*}
    c_K(g) \leq \sum_{\substack{1 \leq d < g, \textrm{ or} \\ d < 2g, \, d \textrm{ even}}}
    \frac{n}{[K(\zeta_d) : \mathbb{Q}(\zeta_d)]}
    + 1 + n \pi(2g)
\end{equation*}
which yields our result. \qed \\

Whilst the above construction does not necessarily improve upon the result given in Theorem \ref{thm:upperDickson} for all fields $K$, one can obtain the following two corollaries: \\

\refstepcounter{mainnum} \label{thm:upperCor1}
\textbf{Corollary \themainnum.} Let $K$ be a primitive abelian number field of (necessarily prime) degree $n$ and conductor $\mathfrak{f}_K$. Then assuming Schinzel's hypothesis H for $K$, we have
\begin{equation*}
    c_K(g) \leq \begin{cases}
    \frac{3}{2} g \big(1 + \frac{n-1}{\mathfrak{f}_K} \big) + 1 + n \pi(2g) & \textrm{ if $\mathfrak{f}_K$ odd, } \\[2mm]
    \frac{3}{2} g \big( 1 + \frac{4(n-1)}{3 \mathfrak{f}_K} \big) + 1 + n \pi(2g) & \textrm{ if $\mathfrak{f}_K$ even. } 
    \end{cases}
\end{equation*}

\textit{Proof.} Since $K$ has conductor $\mathfrak{f}_K$, this implies $[K(\zeta_d) : \mathbb{Q}(\zeta_d)] = 1$ if $\mathfrak{f}_K$ divides $d$, and  $[K(\zeta_d) : \mathbb{Q}(\zeta_d)] = n$ otherwise, noting that $K$ is primitive.  We can therefore easily evaluate the bound given in (\ref{eq:uppermaingenbound}) as
\begin{equation*}
    c_K(g) \leq 
    \sum_{\substack{1 \leq d < g, \textrm{ or} \\ d < 2g, \, d \textrm{ even} \\ \mathfrak{f}_K | d}} \! n \; + 
    \sum_{\substack{1 \leq d < g, \textrm{ or} \\ d < 2g, \, d \textrm{ even} \\ \mathfrak{f}_K \not | d}} \! 1 
    \quad
    + 1 + n \pi(2g)
\end{equation*}
If $\mathfrak{f}_K$ is odd, this evaluates to
\begin{equation*}
    c_K(g) \leq \frac{3}{2} \frac{gn}{\mathfrak{f}_K} + \big( \frac{3}{2} g - \frac{3}{2} \frac{g}{\mathfrak{f}_K} \big) + 1 + n \pi(2g)
\end{equation*}
whilst if $\mathfrak{f}_K$ is even, this yields
\begin{equation*}
    c_K(g) \leq 
    \frac{2gn}{\mathfrak{f}_K}
    + \big( \frac{3}{2} g - \frac{2g}{\mathfrak{f}_K} \big) + 1 + n \pi(2g)
\end{equation*}
which proves the result. \qed \\

\refstepcounter{mainnum} \label{thm:upperCor2}
\textbf{Corollary \themainnum.} Let $K$ be a number field such that its maximal abelian subfield is $\mathbb{Q}$. Then assuming Schinzel’s hypothesis H for $K$, we have $c_K(g) \leq \frac{3}{2} g + n \pi (2g)$. \\

\textit{Proof.}  The above condition implies that $K \cap \mathbb{Q}(\zeta_d) = \mathbb{Q}$ for all $d$, and thus the bound given in (\ref{eq:uppermaingenbound}) implies our result. \qed \\

Finally, it's worth mentioning that the bound given in (\ref{eq:uppermaingenbound}) does not necessarily represent the optimal conditional bound for all genera $g$, even over $\mathbb{Q}$.  For example, under the assumption of Schinzel's hypothesis H, there exist infinitely many integers $k$ such that $k$, $\alpha k - 1$, $\alpha k + 1$, $(\alpha k)^2 + 1$, $(\alpha k)^2 - 2(\alpha k) - 1$, and $(\alpha k)^2 - 2(\alpha k) + 1$ are all prime, where $\alpha := 7!$. \\

From this, one can therefore conditionally construct infinitely many genus $5$ curves $C/\mathbb{Q}$ of the form: 
\begin{align*}
    C_k : y^2 &= x \cdot
    \big(x - (\alpha k)^2(\alpha k-1)(\alpha k+1)\big) \cdot \big(x + (\alpha k)^2(\alpha k-1)(\alpha k+1) \big) \\
    &\qquad \cdot
    \big(x - \alpha k(\alpha k-1)(\alpha k+1)\big) \cdot \big(x + \alpha k(\alpha k-1)(\alpha k+1) \big) \\
    &\qquad \cdot \big(x - (\alpha k-1)(\alpha k+1)\big) \cdot\big(x + (\alpha k-1)(\alpha k+1) \big) \\
    &\qquad \cdot
    \big(x - \alpha k(\alpha k-1)^2\big) \cdot \big(x + \alpha k(\alpha k-1)^2 \big) \\
    &\qquad \cdot
    \big(x - \alpha k(\alpha k+1)^2\big) \cdot\big(x + \alpha k(\alpha k+1)^2 \big) 
\end{align*}
which yields a conditional bound of $c_\mathbb{Q}(5) \leq 10$, and thus one better than the bound of 11 given by (\ref{eq:uppermaingenbound}). \\

Besides the above example, we should mention however that we haven't found any better examples for higher genera over $\mathbb{Q}$, noting that a naive computational search quickly becomes unmanageable for large genus hyperelliptic curves. \\

Finally, I would like to give my sincere thanks to my supervisor, Samir Siksek, for his amazing support and many insightful comments.

\smallskip \footnotesize
\textsc{Mathematics Institute, University of Warwick, CV4 7AL, United Kingdom} \par \nopagebreak
\textit{E-mail address}: \texttt{robin.visser@warwick.ac.uk}
\end{document}